\documentclass[graybox, natbib, nosecnum]{svmult}
\bibpunct{(}{)}{;}{a}{}{,} 

\usepackage{mathptmx}       
\usepackage{helvet}         
\usepackage{courier}        
\usepackage{type1cm}        

\usepackage{makeidx}         
\usepackage{graphicx}        
\usepackage{multicol}        
\usepackage{multirow}
\usepackage[table]{xcolor}
\usepackage[bottom]{footmisc}
\usepackage[normalem]{ulem}	
\usepackage{soul}   
\usepackage{kotex}
\usepackage{amsmath}
\usepackage{amsfonts}
\makeatletter
\newcommand{\tpmod}[1]{\!\!\!\!{\@displayfalse\pmod{#1}}}
\makeatother
\usepackage{mathtools}

\DeclarePairedDelimiter\floor{\lfloor}{\rfloor}
\usepackage[colorlinks=true, allcolors=blue, breaklinks=true]{hyperref}

\usepackage{tikz} 
\usetikzlibrary{positioning,shapes.geometric,arrows,calc}

\def\firstellip{ (1.6, 0)   ellipse [x radius=3cm, y radius=1.5cm, rotate=50]}
\def\secondellip{(0.3, 1cm) ellipse [x radius=3cm, y radius=1.5cm, rotate=50]} 
\def\thirdellip{(-1.6, 0)   ellipse [x radius=3cm, y radius=1.5cm, rotate=-50]} 
\def\fourthellip{(-0.3, 1cm)ellipse [x radius=3cm, y radius=1.5cm, rotate=-50]} 

\usepackage[draft]{changes} 
\definechangesauthor[name={François}, color=orange!50!red]{F}
\definechangesauthor[name={Natan}, color=green!60!blue]{N}

\makeindex             

\begin{document}

\title*{An ethnoarithmetic excursion into the Javanese~calendar}
\author{Natanael Karjanto and Fran\c cois Beauducel}
\institute{N. Karjanto \at Department of Mathematics, University College, Sungkyunkwan University, Natural Science Campus, 2066 Seobu-ro, Suwon 16419, Republic of Korea \hfill \email{natanael@skku.edu}
\and F. Beauducel \at Université de Paris, Institut de physique du globe de Paris, CNRS, 75005 Paris, France \\ 
Institut de recherche pour le développement, Research and Development Technology Center for Geological Disaster, Balai Penyelidikan dan Pengembangan Teknologi Kebencanaan Geologi (BPPTKG), Jl Cendana No. 15, Yogyakarta 55166, Indonesia \hfill \email{beauducel@ipgp.fr}}
%
%
\maketitle
\abstract{A perpetual calendar, a calendar designed to find out the day of the week for a given date, employs a rich arithmetical calculation using congruence. Zeller's congruence is a well-known algorithm to calculate the day of the week for any Julian or Gregorian calendar date. Another rather infamous perpetual calendar has been used for nearly four centuries among Javanese people in Indonesia. This Javanese calendar combines the \emph{Saka} Hindu, lunar Islamic, and western Gregorian calendars. In addition to the regular seven-day, lunar month, and lunar year cycles, it also contains five-day \emph{pasaran}, 35-day \emph{wetonan}, 210-day \emph{pawukon}, octo-year \emph{windu}, and 120-year \emph{kurup} cycles. The Javanese calendar is used for cultural and spiritual purposes, including a decision to tie the knot among couples. In this chapter, we will explore the relationship between mathematics and the culture of Javanese people and how they use their calendar and the arithmetic aspect of it in their daily lives. We also propose an unprecedented congruence formula to compute the \emph{pasaran} day. We hope that this excursion provides an insightful idea that can be adopted for teaching and learning of congruence in number theory.}

\section{Introduction}

Arithmetic and number theory find applications in various cultures throughout the world. In addition to solving everyday problems using elementary arithmetic operations, our ancient ancestors also developed and invented perpetual calendars without any aid of modern electronic calculators and the computer. A perpetual calendar is a system dealing with periods of time that occur repeatedly. It is organized into days, weeks, months, and years. A date of a calendar structure refers to a particular day within that system. A calendar is used for various civil, administrative, commercial, social, and religious purposes. 

The English word ``calendar'' is derived from the Latin word which refers to the first day of the month in the Roman calendar. It is related to the verb \emph{calare} (``to announce solemnly, to call out''), which refers to the ``calling'' of the new moon when it was visible for the first time~\citep{brown1993new}. Another source mentions that the modern English ``calendar'' comes from the Middle English \emph{calender}, which was adopted from the Old French \emph{calendier}. It originated from the Latin word \emph{calendarium}, which meant a ``debt book, account book, register''.  In Ancient Rome, interests were tracked in such books, accounts were settled, and debts were collected on the first day (\emph{calends}) of each month~\citep{stakhov2009mathematics}.
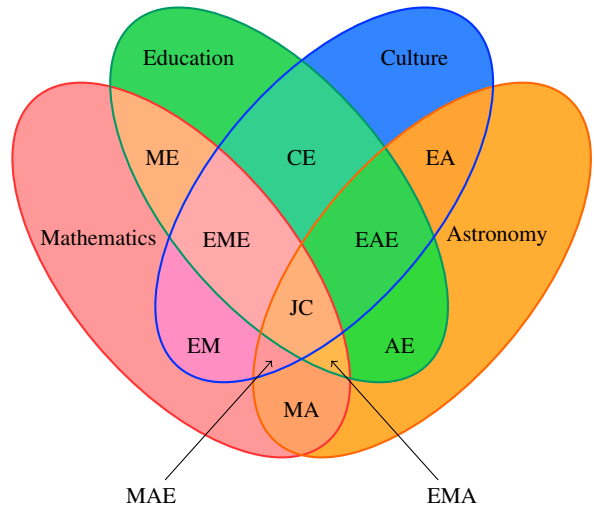
\begin{figure}[h] 
\begin{center}
\begin{tikzpicture}
\begin{scope}[blend group=soft light, opacity = 0.8]
\fill[red!50!white] \thirdellip; 	
\fill[green!80!blue] \fourthellip; 	
\fill[red!40!yellow] \firstellip; 	
\fill[blue!60!cyan] \secondellip;	
\end{scope}
\draw[draw = red!80!white, thick] \thirdellip node [label={[xshift=-1.1cm, yshift=0.1cm] Mathematics}] {};	
\draw[draw = green!60!blue, thick] \fourthellip node [label={[xshift=-1.2cm, yshift=1.5cm] Education}] {};
\draw[draw = red!60!yellow, thick] \firstellip node [label={[xshift=1.0cm, yshift=0.1cm] Astronomy}] {};
\draw[draw = blue!80!cyan, thick] \secondellip node [label={[xshift=1.2cm, yshift=1.5cm] Culture}] {};
		
\node at (-1.85cm,1.5cm){ME};
\node at (-1.3cm,-1cm){EM};
\node at (0.0cm,-1.85cm){MA};
\node at (1.85cm,1.5cm){EA};
\node at (1.3cm,-1cm){AE};
\node at (0.0cm,1.5cm){CE};
\node at (-1.0cm,0.4cm){EME};
\node at (1.0cm,0.4cm){EAE};
\node (MAEX) at (-0.3cm,-1.1cm){};
\node (MAE) at (-2.0cm,-3.0cm){MAE};
\node (EMAY) at ( 0.3cm,-1.1cm){};
\node (EMA) at (2.0cm,-3.0cm){EMA};
\draw[->] (MAE) to (MAEX);
\draw[->] (EMA) to (EMAY);
\node at ( 0.0cm,-0.5cm){JC};
\end{tikzpicture}
\end{center}	
\caption{A theoretical framework for an ethnoarithmetic study of the Javanese calendar. At the first level of the intersection, ME means Mathematics Education, EM means Ethnomathematics, MA means Mathematics Astronomy, AE means Astronomy Education, EA means Ethnoastronomy, and CE means Cultural Education. At the second level of the intersection, EME means Ethnomathematics Education, MAE means Mathematics Astronomy Education, EMA means Ethnomathematics Astronomy, and EAE means Ethnoastronomy Education. The third level of the intersection is the heart of discussion, we have JC as the Javanese calendar.} \label{framework}
\end{figure}

Generally, the calculation of and periods in the calendrical system are synchronized with the cycle of the Sun or the Moon. Hence, the names solar, lunar, and lunisolar calendars, where the latter combined both the Moon phase and the tropical (solar) year. Our current Gregorian and the previously Julian calendars are solar-based calendars, as well as the ancient Egyptian calendar. An example of the lunar calendar is the Islamic calendar. Prominent examples of the lunisolar calendar are the Hebrew, Chinese, Hindu, and Buddhist calendars~\citep{richards1998mapping,dershowitz2008calendrical}.

Studying calendar systems from any culture cannot be separated from attempting to gain insight and understanding not only the arithmetic behind it but also its interconnectivity with mathematics in general, history of mathematics, astronomy, mathematics education, and the sociocultural aspect itself. Hence, the study of perpetual calendars is encompassed by the fields of ethnomathematics and ethnoastronomy, the research areas where diverse cultural groups embrace, practice, and develop mathematics and astronomy into their daily lives. 

The theoretical framework for a discussion on the Javanese calendar considered in this chapter can be summarized in Figure~\ref{framework}. Although the title of this chapter contains the term ``ethnoarithmetic'', indeed that covering a topic on the Javanese calendar encompasses other aspects beyond the cultural and mathematical, in this case, arithmetical aspect per se. As we observe at the third level of the Venn diagram presented in Figure~\ref{framework}, an excursion into the Javanese calendar intersects and includes four distinct but related disciplines: Mathematics, Education, Culture, and Astronomy.

Fundamental work in the area of ethnomathematics has been documented by~\cite{ascher2018mathematics}. The author elaborated on several mathematical ideas and their cultural embedding of people in traditional or small-scale cultures, with some emphasis on time structure and the logic of divination. In particular, the book also covered an explanation of the Balinese calendar. For discussion of the Balinese calendrical system and its multiple cycles, consult~\cite{vickers1990balinese,darling2004marking,proudfoot2007search,ginaya2018balinese}, and~\cite{gislen2019calendars4}.

To the best of our knowledge, the only monograph that discusses an intersection of three components of the theoretical framework is a volume edited by~\cite{rosa2017ethnomathematics}. The book covers diverse approaches and perspectives of ethnomathematics to mathematics education, particularly in several non-Western cultures. Since these topics stimulate debates not only on the nature of mathematical knowledge and the knowledge of a specific cultural group but also on the pedagogy of mathematics classroom, the field of ethnomathematics certainly offers a possibility for improving mathematics education across cultures.

A collection of essays dealing with the mathematical knowledge and beliefs of cultures outside the Western world is compiled by~\cite{selin2000mathematics}. The essays address the connections between mathematics and culture, relate mathematical practices in various cultures, and discuss how mathematical knowledge transferred from East to West. A coverage of calendars in various cultures is also briefly touched.

Since many calendrical systems are invented based on the movement of celestial bodies, particularly the Sun and Moon, we cannot dismantle the role of Astronomy in the study of calendars, including the Javanese one. Various civilizations incorporate the cyclical movements of both the Sun and Moon into their calendrical systems~\citep{ruggles2015calendars}. Some examples of this calendrical-astronomical relationship can be observed amongst others in the Jewish~\citep{cohn2007mathematics}, Indian~\citep{dershowitz2009indian}, Chinese~\citep{martzloff2000chinese,martzloff2016astronomy}, mainland Southeast Asia~\citep{eade1995calendrical}, and Islamic calendars~\citep{proudfoot2006old}. 
	
An explanation of the mathematical and astronomical details of how many calendars function has been covered extensively by~\cite{dershowitz2008calendrical}. A related series with non-Western ethnomathematics is a book on non-Western ethnoastronomy that has been edited by~\cite{selin2000astronomy}. In particular, the book also dedicated one chapter on an astronomical feature of the Javanese calendar, where a season keeper (\emph{pranotomongso} or \emph{pranata mangsa}) guides agricultural activities among rural peasants in Java~\citep{daldjoeni1984pranatamangsa,ammarell1988sky,hidayat2000indo}. The readers might be interested to compare this with the cultural production of Indonesian skylore across three ethnic groups: Banjar Muslim, Meratus Dayak, and Javanese peoples~\citep{ammarell2015cultural}.

Integrating the Javanese calendar into elementary school education as an ethnomathematics study has been attempted by~\cite{utami2020ethnomathematics}.
Another attempt is to embed some topics related to calendrical systems into the first-year seminar on the mathematics of the pre-Columbian Americas~\citep{catepillan2016ethnomathematics}. A study from Taipei, Taiwan, on the movements of the Moon at the primary level, introduced pupils to both the Gregorian and Chinese calendars in the context of physical classroom environments~\citep{hubber2017physical}. 

Although we propose an intersection of only four disciplines in our theoretical framework, the list is by no means exhaustive. Another possibility is to include psychology or the interaction between mathematics and psychology. For instance, a theoretical foundation on the calendrical system and the psychology of time has been modeled by~\cite{rudolph2006fullness}. In particular, the author proposed that Balinese (also equally applied to Javanese) time might be neither `circular' nor `linear', but \emph{profinite}. Since the Javanese calendrical systems involve five-day \emph{pasaran}, seven-day \emph{saptawara/dinapitu}, and 30 \emph{wuku} cycles, it could be modeled by the group of \emph{profinite integers} $\widehat{\mathbb{Z}}$. This group bundles together all different sets of the ring of $p$-adic integers $\mathbb{Z}_p$ and various finite sets $\mathbb{Z}/(p)$ of integers modulo $p$, where $p$ denotes prime integers~\citep{milneCFT}.

A discussion on calendar timekeeping scheme in Southeast Asia is covered by~\cite{gislen2018lunisolar}. An overview of the calendars in Southeast Asia is given by~\cite{gislen2019calendars1}. In their subsequent papers, they and L{\^a}n also discussed calendars from Burma, Thailand, Laos, and Cambodia~\citep{gislen2019calendars2}, Vietnam~\citep{lan2019calendars}, Malaysia and Indonesia~\citep{gislen2019calendars4}, eclipse calculation~\citep{gislen2019calendars5}, and chronicle inscriptions~\citep{gislen2019calendars6}. See also~\cite{eade1995calendrical,ohashi2009mainland}, and~\citep{golzio2012calendar} for basic facts and further explanations of the calendrical systems in India and (Mainland) Southeast Asia.

This chapter is organized as follows. After this Introduction, the following section briefly covers the calendars from ancient and contemporary times, which includes the pre-Gregorian and Gregorian calendars. After that, a discussion on the Javanese calendar in the context of ethnoarithmetic will be dedicated exclusively in one section. The final section concludes our discussion.

\section{Ancient and modern calendars}			\label{ancientmodern}

This section briefly covers the calendar in the pre-Gregorian era and Zeller's congruence algorithm in the Gregorian calendar.

\subsection{Pre-Gregorian calendars}

Before the Gregorian calendar that we are using today, numerous calendars have been used in various parts of the world. The Egyptian calendar was among the first solar calendars with its history dated back to the fourth millennium~BCE.

A Mesoamerican civilization developed by the Maya peoples also developed a calendar system where the year was divided into 18 months of 20 days~\citep{ascher2018mathematics,stakhov2009mathematics}. The Mayan calendar system contains three separate calendars, and although they are not related mathematically, all of them are linked in a single calendar system. The first one is called the long count, used by the Mayans to measure date chronology for history recording. The second one is the \emph{Haab} calendar, a non-chronological civil calendar. The third one is the Mayan religious calendar, also non-chronological, called the \emph{Tzolkin}~\citep{ascher2018mathematics,cohn2007mathematics}.

Related to the Maya civilization is the Inca Empire, where the latter also developed its calendar system, which was based on several different astronomical cycles, including the solar year, the synodic and sidereal lunar cycles, and the local zenith period~\citep{urton2010social}. A variety of mathematical development among the native Americans from the prehistoric to present has been compiled by~\cite{closs1996native}. In particular, one chapter of the book discusses the calendrical system of the Nuu-chah-nulth (formerly referred to as the Nootka), one of the indigenous peoples of the Pacific Northwest Coast in Canada~\citep{folan1986calendrical}.

Some calendars are based on lunisolar, and one of them is the Jewish calendar. Also called the Hebrew calendar, it is still used until today, mainly for Jewish religious observance. Although the Jewish calendar was developed in its current format during the Talmud period in the fifth century CE by Rabbi Hillel II, the up to date counting has accumulated more than 5000 years~\citep{ascher2018mathematics,cohn2007mathematics}. For our information, AM 5781 will begin at sunset on Friday, 18 September 2020, and will end at sunset on Monday, 6 September 2021~CE. Here, AM means \emph{Anno Mundi}, the Latin phrase for ``in the year of the world''.

Although the traditional Chinese calendar is also a lunisolar type, the most fundamental component is the sexagenary cycle, a cycle of sixty terms marked by coordination between 12 celestial stems and ten terrestrial branches, which is also known as the Stems-and-Branches or \emph{g\={a}nzh\={\i}} (干支)~\citep{martzloff2000chinese}. A discussion on the mathematical aspect of the Chinese calendar has been covered extensively by~\cite{aslaksen2001fake,aslaksen2009when,aslaksen2010mathematics}. A historical aspect of the Chinese calendar is discussed by~\cite{sun2015chinese}. Astronomical aspects and the mathematical structures underlying the calculation techniques of the Chinese calendar are highlighted by~\cite{martzloff2016astronomy}. 

The \emph{Hijri}, or Islamic calendar, is a lunar calendar consisting of 12 lunar (synodic) months in a year of 354 or 355 days. It is still widely used in predominantly Muslim countries alongside the Gregorian calendar, primarily for religious purposes. The current counting started from 622 CE, commemorating the emigration of Prophet Muhammad and his companions from Mecca to Medina (\emph{hijra})~\citep{hassan2017muslimcalendar}. In the Gregorian calendar reckoning, the current Islamic year is 1422 AH, which approximately runs from 20 August 2020 until 9 August 2021. AH means \emph{Anno Hegirae}, the Latin phrase for ``in the year of the Hijra''. In particular, the history of Muslim calendars in Southeast Asia in the historical and cultural context has been successfully traced by~\cite{proudfoot2006old}.

\subsection{Gregorian calendar}

The language of congruences is a basic building block in arithmetic and number theory. It allows us to operate with divisibility relationships in a similar way as we deal with equalities. Congruences have many applications, and one of them is to determine the day of the week for any date for a given perpetual calendar. In particular, the procedure for finding the day of the week for a given date in the Gregorian calendar has been discussed by~\cite{carroll1887find,conway1973tomorrow,rosen2011elementary}. See also~\cite{gardner2007universe,cohen2000day}.

The Gregorian calendar was promulgated by Pope Gregory XIII and is a replacement of the Julian calendar, proposed by Julius Caesar in 46~BC. Starting from 1582, the Catholic states were among the first countries to adopt the Gregorian calendar by skipping ten days in October, Thursday, 4 October 1582 was followed by Friday, 15 October 1582. Greece and Turkey are among the last countries that adopted the Gregorian calendar when they changed on 1 March 1923 and 1 January 1926, respectively.

Let $W$ denote the day of the week from Saturday = $0$ to Friday = $6$, $k$ denote the day of the month, $m$ denote month, and $N$ denote year. For the month, the convention is March = $3$, April = $4$, \dots, but January = $13$, and February = $14$. For the year, $N$ is the current year unless the month is January or February, for which $N$ is the previous year. The relationship between year and century is given by $N = 100 C + Y$, where $C$ denotes zero-based century and $Y$ denotes a two-digit year. Note that the purpose of adopting this expression is not to confuse $C$ with the standard century number, which is $C + 1$ for the first 99 years. The formula for finding the day of the week $W$ of day $k$ of month $m$ of year $N$ is given by Zeller's congruence algorithm~\citep{zeller1882grundaufgaben,zeller1883problema,zeller1885kalender,zeller1887kalender}:
\begin{equation}
W \equiv \left(k + \floor*{\frac{13(m + 1)}{5}} + Y + \floor*{\frac{Y}{4}} + \floor*{\frac{C}{4}} - 2C \right) \tpmod 7.
\label{weekgreg}
\end{equation}
Here, $\floor{x}$ denotes the floor function or greatest integer function and mod is the modulo operation or remainder after division.

As an example, we are interested in finding the days of the week when a Javanese heroine and educator Raden Adjeng Kartini was born and passed away. She was born in Jepara, a town on the north coast of Java, around 80~km north-east direction from Semarang, the present capital of Central Java province in Indonesia. Kartini was a pioneer for girls' education and women's emancipation rights in Indonesia, at the time when Indonesia was still a part of the Dutch East Indies colonial empire. She could be correlated with her European counterparts, including an English advocate of women's right Mary Wollstonecraft (27 April 1759--10 September 1797) or a Finnish social activist Minna Canth (19 March 1844--12 May 1897)~\citep[cf.][]{kartini1911door,kartini1985letters,wollstonecraft1995wollstonecraft,canth2013tyomiehen}. 

Kartini was born on 21 April 1879, so we have $k = 21$, $m = 4$, $N = 1879$, $C = 18$, and $Y = 79$. She passed away on 17 September 1904, at the age of 25, and we identify that $k = 17$, $m = 9$, $N = 1904$, $C = 19$, and $Y = 4$. Using Zeller's congruence formula~\eqref{weekgreg}, the days when Kartini was born $W_b$ and passed away $W_d$ can be calculated as follows, respectively:
\begin{align*}
W_b &\equiv \left(21 + \floor*{\frac{13 (5)}{5}} + 79 + \floor*{\frac{79}{4}} + \floor*{\frac{18}{4}} - 36 \right) \tpmod {7}  \\
    &\equiv \left(21 + 13 + 79 + 19 + 4 - 36 \right) \tpmod{7} \equiv 100 \, \tpmod{7} \equiv 2 \, \tpmod {7} \\
W_d &\equiv \left(17 + \floor*{\frac{13 (10)}{5}} + 4 + \floor*{\frac{4}{4}} + \floor*{\frac{19}{4}} - 38 \right) \tpmod{7} \\
    &\equiv \left(17 + 26 + 4 + 1 + 4 - 38 \right) \tpmod{7} \equiv 14 \, \tpmod{7} \equiv 0 \, \tpmod{7}.
\end{align*}
Hence, Kartini was born on Monday, 21 April 1879 and passed away on Saturday, 17 September 1904. She was buried at Bulu Village, Rembang, Central Java, around 100~km east of Jepara.

\section{Javanese calendar}
\label{javanese}

This section features the main characteristics of the Javanese calendar. We start with the geographical location of the island of Java, the Javanese people, and a brief historical background of the calendar. After providing a detailed discussion, we close the section by computer implementation of the Javanese calendar.

\subsection{Where is Java?}

\begin{figure}[h] 
\includegraphics[trim=20 180 20 80,clip,width=\textwidth]{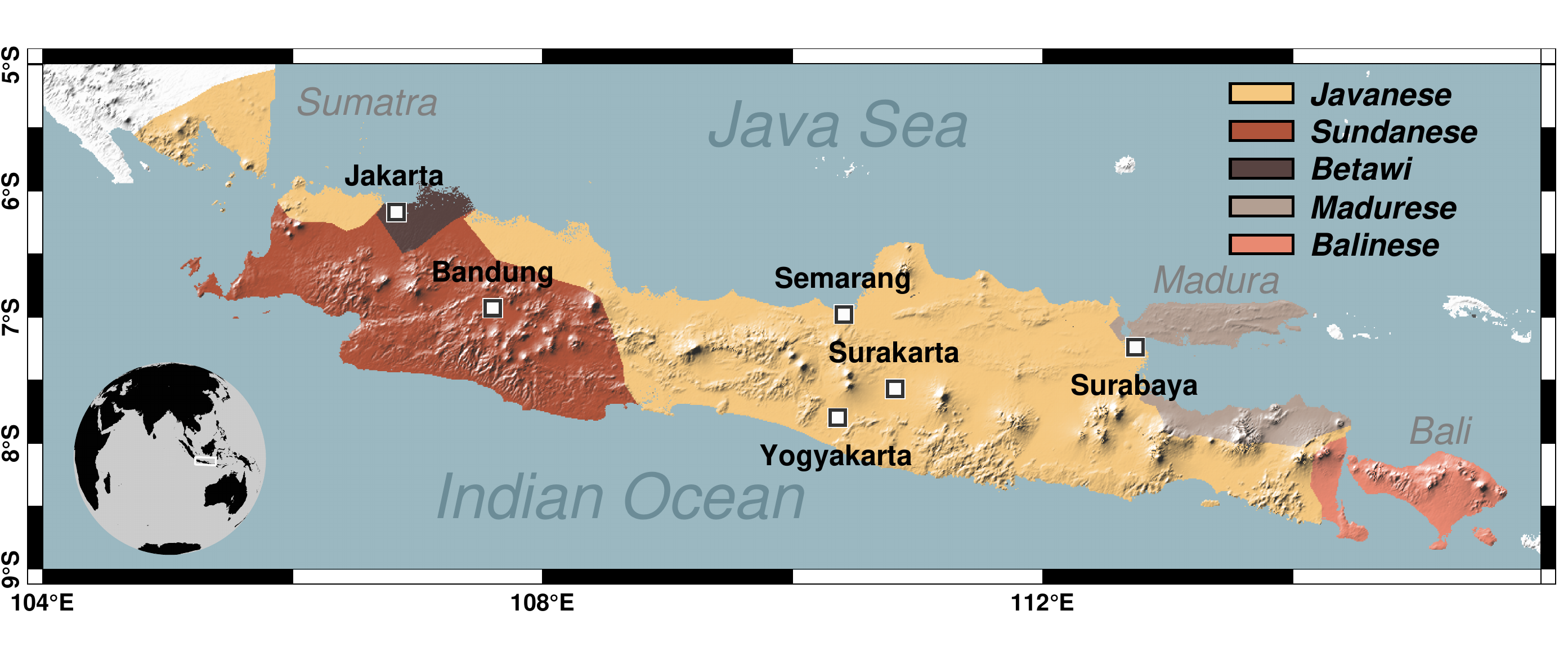}
\caption{The situation of the Java island, main cities, and present spoken languages: Javanese (Central and East Java, and a small enclave in North-West Java), Sundanese (West Java), Betawi (in and around Jakarta metropolitan area), Madurese (Madura Island and a part of North-Eastern Java, and Balinese (Bali Island and a small part of Eastern Java). Basemap uses ETOPO5 and SRTM3 topographic data and shaded relief mapping code~\citep{BeauducelDEM20}. ETOPO5 is a five arc-minute resolution relief model for the Earth's surface that integrates land topography and ocean bathymetry dataset. SRTM3, the Shuttle Radar Topography Mission, is a three arc-second resolution digital topographic database of land elevation limited to latitudes from 60$^\circ$ South to 60$^\circ$ North.}	\label{javamap}
\end{figure}

Java is an island in Indonesia, not a programming language, although the latter was renamed after some Javanese coffee by its founders. It is located in the southern hemisphere, around 800~km from the Equator. It extends from latitude 6$^\circ$ to 8$^\circ$ South, and longitude 105$^\circ$ to 114$^\circ$ East. With an area of 150,000 square kilometers, it is about 1000~km long from west to east and around 200~km wide from north to south. The island lies between Sumatra to the west and Bali to the east. It is bordered by the Java Sea on the north and the Indian Ocean on the south (see Figure \ref{javamap}). It is the world's 13$^\text{th}$ largest island~\citep{cribb2013historical}.

\subsection{Who are Javanese people?}

The Javanese people are a native ethnic group to the island of Java. They form the largest ethnic group in Indonesia, with more than 95 million people live in Indonesia and approximately 5 million people live abroad. Although they predominantly reside in the central and eastern parts of the island, they are also scattered in various parts of the country~\citep{taylor2003indonesia,ananta2015demography}. The Javanese people possess and speak a distinct language from Indonesian, called the Javanese language, a member of the Austronesian family of languages written in Javanese script \emph{hanacaraka} or \emph{dentawyanjana}. Thanks to its long history and legacy of Hinduism and Buddhism in Java, the language adopted a large number of Sanskrit words~\citep{marr1986southeast,errington1998shiting}.

\subsection{A background of the Javanese calendar}

Javanese people use the Javanese calendar simultaneously with two other perpetual calendars, the Gregorian and Islamic calendars. The former is the official calendar of the Republic of Indonesia and the latter is used mainly for religious purposes. Prior to the adoption of the Javanese calendar in 1633~CE, Javanese people used a calendrical system based on the lunisolar Hindu \emph{Saka} calendar~\citep{ricklefs1993history}. 
\begin{figure} 
\includegraphics[width=\textwidth]{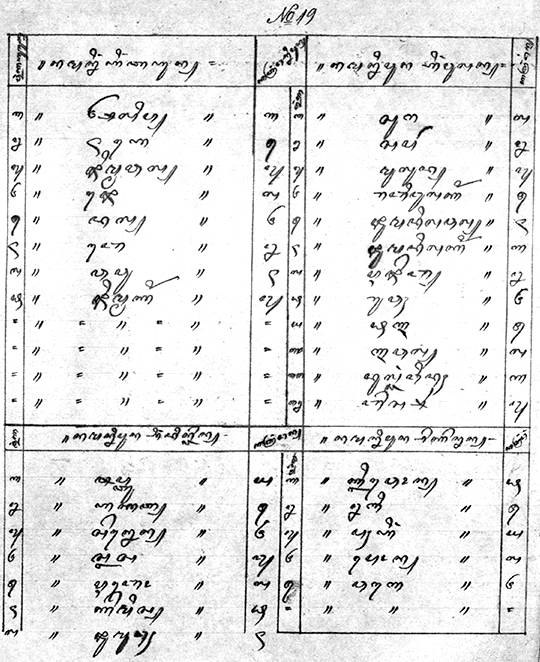}
\caption{Table list of 8 \emph{taun}, 12 \emph{wulan}, 7 \emph{dinapitu}, and 5 \emph{pasaran} names in Javanese Sanskrit, after \cite{warsapradongga1892}.}
\label{hanacaraka}
\end{figure}

The Javanese calendar was inaugurated by Sultan Agung Adi Prabu Hanyakraku-suma~(1593--1645~CE), or simply Sultan Agung, the third Sultan of Mataram who ruled Central Java from 1613~CE until 1645~CE. Although the counting of the year follows the Saka calendar, the Javanese calendar employs a similar lunar year as the Hijri calendar instead of the solar year system like the former~\citep{gislen2019calendars1}. The Javanese calendar is sometimes referred to as AJ (\emph{Anno Javanico}), the Latin phrase for Javanese Year. Since 2008, the difference between the Gregorian and Javanese calendars is about 67 years, where the current year 2020~CE corresponds to 1953~AJ~\citep{oey2001java,raffles2018history}.

\begin{table}
\caption{Main cycles of the Javanese calendar.}
\label{cycles}
\begin{tabular}{p{2cm}p{2.4cm}p{2cm}p{4.9cm}}
\hline\noalign{\smallskip}
Cycle Name & Length & Unit & Comment\\
\noalign{\smallskip}\svhline\noalign{\smallskip}
Pancawara & 5 & day & Javanese week\\
 & & & 5 {\it pasaran} names (see Table \ref{pasaran})\\
\noalign{\smallskip}
Wuku & 7 & day & Gregorian/Islamic week\\
 & & & 7 {\it dinapitu} names (see Table \ref{dinapitu})\\  
 & & & 30 {\it wuku} names\\
\noalign{\smallskip}
{Wetonan} & 35 & day & 35-day names as {\it Dinapitu} and {\it Pasaran}\\
\noalign{\smallskip}
Wulan & 29 or 30 & day & day number in a {\it wulan} is {\it dina}\\
 & & & 12 {\it wulan} names (see Table \ref{wulan})\\
\noalign{\smallskip}
{Pawukon} & 30 & {\it Wuku} & 30 weeks $\equiv$ 210 days\\
\noalign{\smallskip}
Taun & 12 & {\it Wulan} & 354 or 355 days\\
 & & & {\it taun} number starts on 1555~AJ\\
 & & & 8 {\it taun} names (see Table \ref{taun})\\
\noalign{\smallskip}
Windu & 8 & {\it Taun} & 96 {\it wulan} $\equiv$ 81 {\it Wetonan} $\equiv$ 2,835 days\\
 & 4 & {\it Windu} & 32 \emph{taun}, 4 \emph{windu} names\\
\noalign{\smallskip}
{Lambang} & 2 & {\it Windu} & 16 \emph{taun}, 2 \emph{lambang} names\\
\noalign{\smallskip}
Kurup & 15 & {\it Windu} & 120 {\it Taun} $-1$ day $\equiv$ 42,524 days\\
 & & & 4 {\it kurup} names until today (see Table \ref{kurup})\\
\noalign{\smallskip}\hline\noalign{\smallskip}
\end{tabular}
\end{table}

\subsection{Some characteristics of the Javanese calendar}

Different from many other calendars that employ a seven-day week cycle, the Javanese calendar adopts a five-day week cycle, known as \emph{pancawara}. Amalgamating with the seven-day week cycle of the Gregorian and Islamic calendars, namely the \emph{saptawara} cycle, one obtain the 35-day cycle, known as \emph{wetonan}~\citep{darling2004marking}. This foundation cycle interferes with additional cycles:
\begin{itemize}
\item a 210-day cycle of 30 weeks, named as the \emph{pawukon};
\item a more complex combination of the lunar month \emph{wulan}, which has 29 or 30 days, the lunar year \emph{taun}, which is 12 lunar months, \emph{windu}, which is eight lunar years, and finally \emph{kurup} of 15 \emph{windu} or 120 lunar years minus one day, which matches exactly the Islamic calendar cycle~\citep{proudfoot2007search}.
\end{itemize}

There are also additional cycles, but they are no longer used in the Javanese tradition. To be exhaustive, we list them here \citep{richmond1956time,zerubavel1989seven}:
\begin{itemize}
\item a six-day cycle called the \emph{Paringkelan}: `Tungle', `Aryang', `Wurukung', `Paningron', `Uwas', and `Mawulu';
\item an eight-day cycle called the \emph{Padewan}: `Sri', `Indra', `Guru', `Yama', `Rudra', `Brama', `Kala', and `Uma';
\item a nine-day cycle called the \emph{Padangon}: `Dangu', `Jagur', `Gigis', `Kerangan', `Nohan', `Wogan', `Tulus', `Wurung', and `Dadi';
\end{itemize}
The list of main cycles and characteristics are summarized in Table~\ref{cycles} and are detailed in the following subsubsections. Figure~\ref{hanacaraka} shows an original table of the main cycles' names written in the Javanese script \emph{hanacaraka}.

\begin{table}[h] 
\caption{Names of the five \emph{pasaran} days in the Javanese week and commonly associated symbols.}
\label{pasaran}
\centering
\begin{tabular}{p{1.5cm}p{1.5cm}p{1.6cm}p{1.6cm}p{1.6cm}p{1.6cm}p{1.6cm}}
\hline\noalign{\smallskip}
Ngoko & Krama & Meaning & Element & Color & Direction & Posture\\
\noalign{\smallskip}\svhline\noalign{\smallskip}
Pon & Petak & -- & water & yellow & West & sleep\\
Wage & Cemeng & dark & earth & black & North & sit down\\
Kliwon & Asih & affection & spirit & mixed color & focus/centre & stand-up\\
Lêgi & Manis & sweet & air & white & East & turn back\\
Pahing & Pahit & bitter & fire & red & South & to face\\
\noalign{\smallskip}\hline\noalign{\smallskip}
\end{tabular}
\end{table}
\subsubsection{Pancawara and Pasaran}

The Javanese five-day cycle is named \emph{pancawara} and made of five days known as \emph{pasaran}: `Pon', `Wage', `Kliwon', `Lêgi', and `Pahing' (or `Paing'). The word comes from `pasar' which means market. Historically, the market was held and operated on a five-day cycle based on a \emph{pasaran} day, e.g., `Pasar Kliwon', or `Pasar Legi'. Until today, most of the markets in Java still have a \emph{pasaran} name like the Kliwon Market in Kudus, Central Java, although they usually operate every day~\citep{oey2001java}.

Each \emph{pasaran} is associated with some symbols, in particular, the five classical elements of Aristotle, colors, cardinal directions (note that the Javanese culture has five of them including a Center), and human posture~\citep{pigeaud1977javanese},~\citep[cf.][]{brinton1893native}. See Table~\ref{pasaran}.

The computation of the \emph{pasaran} day from the Gregorian or Islamic calendar date can be performed using a similar strategy as Zeller's congruence algorithm. Indeed, the formula is based on the observation that the day of the week progresses predictably based upon each subpart of the date, i.e., day, month, and year. In this case, we will consider the five-day Javanese week \emph{pasaran}. Each term within the formula is used to calculate the offset needed to obtain the correct day.

Let $P$ denote the \emph{pasaran} day, in the order of Pon~$= 0$, Wage~$= 1$, Kliwon~$= 2$, Lêgi~$= 3$, and Pahing~$= 4$. We now propose the following \emph{pasaran} day congruence formula:
\begin{equation}
P \equiv \left(k + \floor*{\frac{3(m + 1)}{5}} + \floor*{\frac{Y}{4}}+ 4C - 4\floor*{\frac{C}{4}}\right) \tpmod 5.
\label{weekjav}
\end{equation}

\noindent where $k$, $m$, $C$, and $Y$ denote the same variables as the ones defined for Zeller's congruence formula (see expression~\eqref{weekgreg} in the previous section). Each term in~\eqref{weekjav} can be analyzed as follows:
\begin{itemize}
\item $k$ represents the progression of the day of the week based on the day of the month since each successive day results in an additional offset of one;
\item the term $+\floor*{\frac{3(m + 1)}{5}}$ adjusts for the variation in the days of the month. Indeed, starting from March to February, the days in a month are \{31, 30, 31, 30, 31, 31, 30, 31, 30, 31, 31, 28/29\}. The last element, February's 28/29 days, is not a problem since the formula had rolled it to the end. The number of days for the 11 first elements of this sequence modulo $5$ (still starting with March) would be \{1, 0, 1, 0, 1, 1, 0, 1, 0, 1, 1\} which basically alternates the subsequence \{1, 0, 1, 0, 1\} every five months and give the number of days that should be added to the next month. The fraction $\frac{3}{5} \equiv 0.6$ applied to $m+1$ and the floor function have that effect and will add the proper amount of days;
\item since there are 365 days in a non-leap year, and $365 \tpmod{5} \equiv 0$, there is no need to add an offset for the normal year;
\item for the leap years, $366 \tpmod{5} \equiv 1$ so a day must be added to the offset value by the term $+\left\lfloor {\frac{Y}{4}}\right\rfloor$;
\item there are $36,524$ days in a normal century and $36,525$ days in each century divisible by $400$; the term $4C$ adds $36,524\tpmod{5} \equiv 4$ days for any century, and $-4\floor*{\frac{C}{4}}$ removes these $4$ days for a century divisible by $400$;
\item the overall function, $\tpmod{5}$, normalizes the result to reside in the range from $0$ to $4$, which yields the index of the correct day for the date being analyzed.
\end{itemize}

For example, to find the \emph{pasaran} of the very first day in the Javanese calendar (1 Sura 1555~AJ Alip) which corresponds to 8 July 1633~CE, we have $k = 8$, $m = 7$, $C = 16$, and $Y = 33$. Using~\eqref{weekjav}, we obtain $P \equiv 8 + 4 + 8 + 64 - 16 \equiv 68$ $\tpmod{5}$ $\equiv 3$. Hence, 8 July 1633~CE was a `Lêgi'. From the previously considered example, the \emph{pasaran} days when Kartini was born $P_b$ and passed away $P_d$ can be calculated as follows:
\begin{align*}
P_b \equiv \left(21 + 3 + 19 + 72 - 16 \right) \tpmod 5 \equiv 99 \, \tpmod 5 &\equiv 4 \, \tpmod 5 \\
P_d \equiv \left(17 + 6 +  1 + 76 - 16 \right) \tpmod 5 \equiv 84 \, \tpmod 5 &\equiv 4 \, \tpmod 5.
\end{align*}
Hence, both 21 April 1879 and 17 September 1904 fall on a Pahing.
\begin{figure}[h] 
\begin{center}
\includegraphics[width=0.85\textwidth]{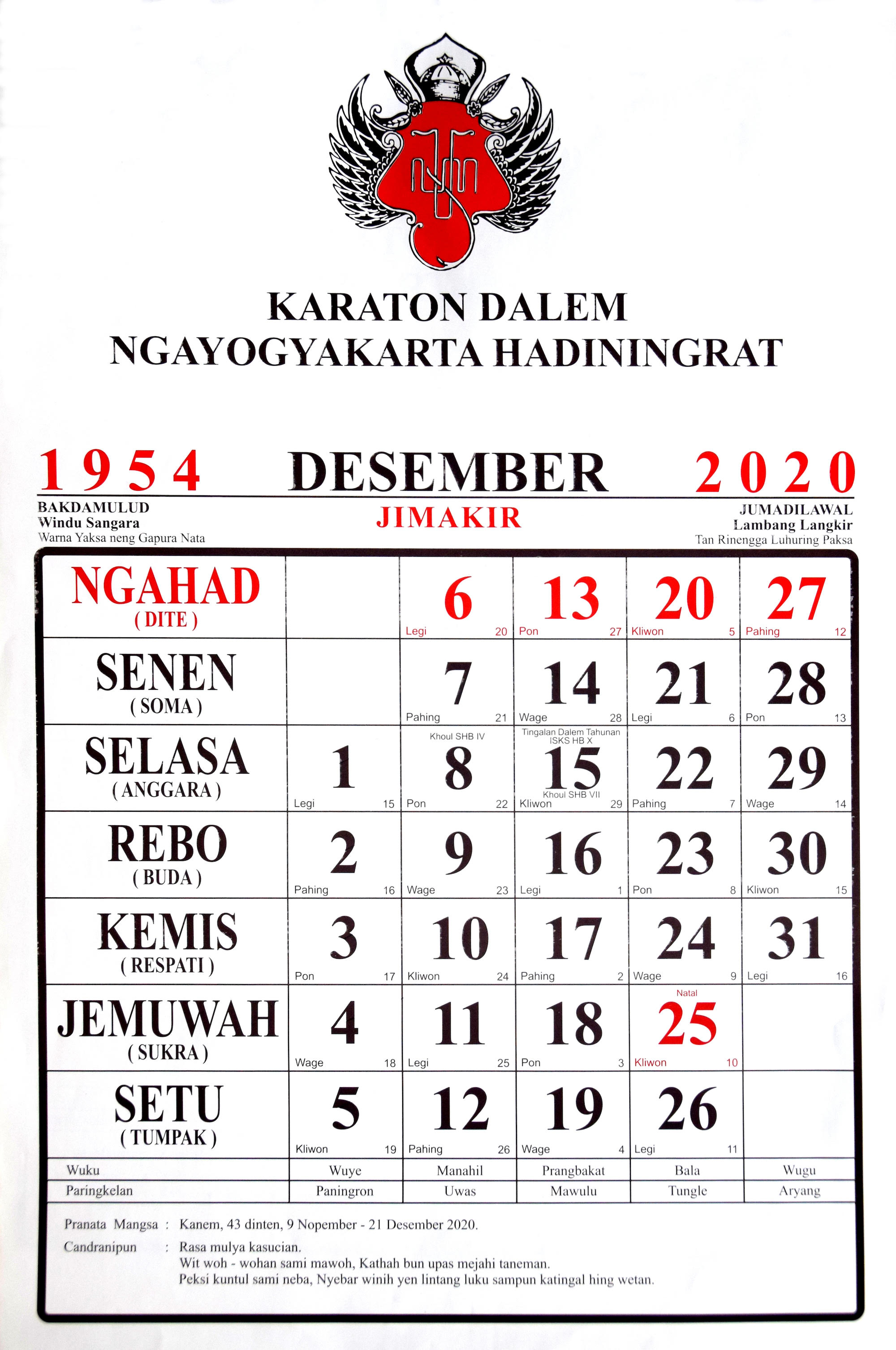}
\end{center}	
\caption{An example of the Javanese calendar for December 2020~CE issued by the Kraton palace in Yogyakarta. It contains information on \emph{Jimakir taun}, \emph{Bakdamulud} and \emph{Jumadilawal wulans}, \emph{pasaran} day, \emph{padinan} weekday, \emph{wuku}, and \emph{paringkelan} amongst others.} \label{2020calendar}
\end{figure}

\subsubsection{Dinapitu, Wuku and Pawukon}
\emph{Dinapitu} literally means `day seven' in Javanese, and corresponds to the day names in the Gregorian/Islamic calendar week with exact equivalence, i.e., from Monday to Sunday: `Sênèn', `Selasa' (or `Slasa'), `Rêbo', `Kêmis', `Jemuwah' (or `Jumungah'), `Sêtu', and  `Ngahad' (or `Ahad'). In recent literature, we often find the name of the days in the Indonesian language (from Monday to Sunday: `Senin', `Selasa', `Rabu', `Kamis', `Jumat', `Sabtu', and `Minggu'). They are also associated with particular symbolic meanings. See Table~\ref{dinapitu}. The names of the days of the week are similar in both languages as they are absorbed from Arabic except for Sunday in Indonesian which was assimilated from the Portuguese `Domingo'.
\begin{table} 
\caption{List of names of the seven days in a \emph{wuku} Gregorian week~\citep[cf.][]{rizzo2020what}.}
\label{dinapitu}
\begin{tabular}{p{2.8cm}p{2.8cm}p{2.8cm}p{2.8cm}}
\hline\noalign{\smallskip}
Dinapitu & Padinan & Week Day & Symbol\\
\noalign{\smallskip}\svhline\noalign{\smallskip}
Ngahad & Dite & Sunday & silent\\
Sênèn & Soma & Monday & forward\\
Selasa & Anggara & Tuesday & backward\\
Rêbo & Buda & Wednesday & turn left\\
Kêmis & Respati & Thursday & turn right\\
Jemuwah & Sukra & Friday & up\\
Sêtu & Tumpak & Saturday & down\\
\noalign{\smallskip}\hline\noalign{\smallskip}
\end{tabular}
\end{table}

The seven-day cycle is called \emph{saptawara} or \emph{padinan} and a week is named a \emph{wuku}. The period of 30 \emph{wuku} makes the \emph{pawukon} cycle, i.e., 210 days~\citep{proudfoot2007search}. There are 30 different names of \emph{wuku}~\citep{soebardi1965calendrical,headley2004javanese}. See Table~\ref{wuku}. Figure~\ref{2020calendar} shows an example of the Javanese calendar with \emph{padinan}, \emph{pasaran}, \emph{wuku}, and \emph{paringkelan} information amongst others.
\begin{table}
\caption{List of names of the 30 different \emph{wuku} in a \emph{pawukon}.}
\label{wuku}
\begin{tabular}{p{0.3cm}p{1.5cm}p{0.3cm}p{1.8cm}p{0.3cm}p{1.8cm}p{0.3cm}p{2.2cm}p{0.3cm}p{1.8cm}}
\hline\noalign{\smallskip}
1 & Sinta   & 7 & Warigalit & 13 & Langkir & 19 & Tambir & 25 & Bala\\
2 & Landep  & 8 & Warigagung & 14 & Mandasiya & 20 & Medangkungan & 26 & Wugu\\
3 & Wukir   & 9 & Julungwangi & 15 & Julungpujut & 21 & Maktal & 27 & Wayang\\
4 & Kurantil & 10 & Sungsang & 16 & Pahang & 22 & Wuye  & 28 & Kulawu\\
5 & Tolu    & 11 & Galungan & 17 & Kuruwelut & 23 & Manahil & 29 & Dukut\\
6 & Gumbreg & 12 & Kuningan & 18 & Marakeh & 24 & Prangbakat & 30 & Watugunung\\
\noalign{\smallskip}\hline\noalign{\smallskip}
\end{tabular}
\end{table}

\subsubsection{Wetonan}

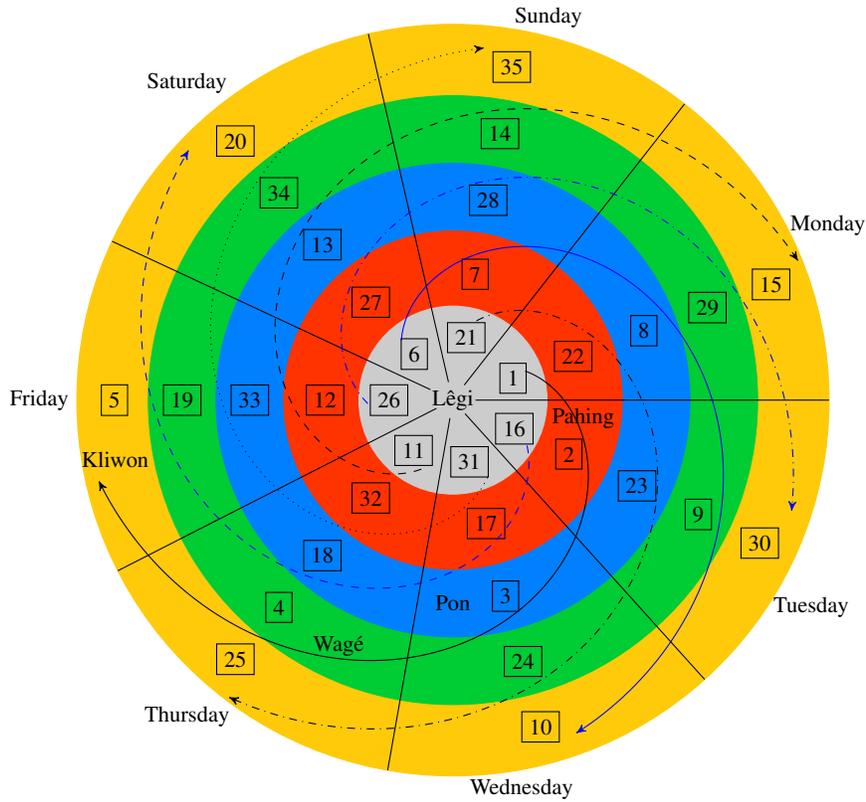
\begin{figure}[h!] 
\begin{center}
\begin{tikzpicture}[auto, main node/.style={draw}]
\def\R{5cm} 
\begin{scope}
\filldraw[yellow!60!orange] (0,0) circle[radius=\R];
\filldraw[green!80!blue] (0,0) circle[radius=4.05*\R/5];
\filldraw[blue!50!cyan] (0,0) circle[radius=3.15*\R/5];
\filldraw[red!80!yellow] (0,0) circle[radius=2.25*\R/5];
\filldraw[black!20!white] (0,0) circle[radius=1.25*\R/5];
\end{scope}
\draw (0.3,0) -- (  0:\R)
(0.15,0.2) -- ( 52:\R)    
(-0.06,0.24) -- (103:\R)     
(-0.30,0.14) -- (155:\R)
(-0.22,-0.11) -- (207:\R)
(-0.04,-0.24) -- (260:\R)
(0.16,-0.18) -- (312:\R);
\node at ( 0: 0.0*\R/5) {L\^{e}gi};
\node at (352: 1.75*\R/5) {Pahing};
\node at (270: 2.7*\R/5) {Pon};
\node at (245: 3.6*\R/5) {Wag\'{e}};
\node at (190: 4.56*\R/5) {Kliwon};
\node at (25: 1.1*\R) {Monday};
\node at (76: 1.05*\R) {Sunday}; 
\node at (130: 1.1*\R) {Saturday}; 
\node at (180: 1.1*\R) {Friday}; 
\node at (230: 1.1*\R) {Thursday};
\node at (280: 1.05*\R) {Wednesday}; 
\node at (330: 1.1*\R) {Tuesday}; 
		
\node[main node] (1) at ( 20: 0.85*\R/5) {$1$};
\node[main node] (2) at (335: 1.7*\R/5) {$2$};
\node[main node] (3) at (285: 2.7*\R/5) {$3$};
\node[main node] (4) at (230: 3.6*\R/5) {$4$};
\node[main node] (5) at (180: 4.5*\R/5) {$5$};
		
\node[main node] (6) at (130: 0.8*\R/5) {$6$};
\node[main node] (7) at ( 80: 1.7*\R/5) {$7$};
\node[main node] (8) at ( 20: 2.7*\R/5) {$8$};
\node[main node] (9) at (335: 3.6*\R/5) {$9$};
\node[main node] (10) at (285: 4.5*\R/5) {$10$};
		
\node[main node] (11) at (232: 0.85*\R/5) {$11$};
\node[main node] (12) at (180: 1.7*\R/5) {$12$};
\node[main node] (13) at (130: 2.7*\R/5) {$13$};
\node[main node] (14) at ( 80: 3.6*\R/5) {$14$};
\node[main node] (15) at ( 20: 4.5*\R/5) {$15$};
		
\node[main node] (16) at (335: 0.9*\R/5) {$16$};
\node[main node] (17) at (285: 1.7*\R/5) {$17$};
\node[main node] (18) at (230: 2.7*\R/5) {$18$};
\node[main node] (19) at (180: 3.6*\R/5) {$19$};
\node[main node] (20) at (130: 4.5*\R/5) {$20$};
	
\node[main node] (21) at (78: 0.85*\R/5) {$21$};
\node[main node] (22) at ( 20: 1.7*\R/5) {$22$};
\node[main node] (23) at (335: 2.7*\R/5) {$23$};
\node[main node] (24) at (285: 3.6*\R/5) {$24$};
\node[main node] (25) at (230: 4.5*\R/5) {$25$};
		
\node[main node] (26) at (180: 0.85*\R/5) {$26$};
\node[main node] (27) at (130: 1.7*\R/5) {$27$};
\node[main node] (28) at ( 80: 2.7*\R/5) {$28$};
\node[main node] (29) at ( 20: 3.6*\R/5) {$29$};
\node[main node] (30) at (335: 4.5*\R/5) {$30$};
		
\node[main node] (31) at (285: 0.85*\R/5) {$31$};
\node[main node] (32) at (230: 1.7*\R/5) {$32$};
\node[main node] (33) at (180: 2.7*\R/5) {$33$};
\node[main node] (34) at (130: 3.6*\R/5) {$34$};
\node[main node] (35) at ( 80: 4.5*\R/5) {$35$};
		
\draw [->, >=stealth', domain= 5.9:9.2,variable=\t,smooth,samples=75]
plot ({-\t r}: {1.15*(\t - 5)});
		
\draw [->, >=stealth', blue, domain= 4:7.5,variable=\t,smooth,samples=75]
plot ({-\t r}: {1.05*(\t - 3)});
		
\draw [->, >=stealth', dashed, domain= 2:5.9,variable=\t,smooth,samples=75]
plot ({-\t r}: {1.01*(\t - 1)});
		
\draw [->, >=stealth', blue, dashed, domain= 0.55:3.9,variable=\t,smooth,samples=75]
plot ({-\t r}: {1.1*(\t + 0.5)});
		
\draw [->, >=stealth', dashdotted, domain= 5:8.5,variable=\t,smooth,samples=75]
plot ({-\t r}: {1.1*(\t - 4)});

\draw [->, >=stealth', blue, dashdotted, domain= 3.1:6.6,variable=\t,smooth,samples=75]
plot ({-\t r}: {1.03*(\t - 2)});
		
\draw [->, >=stealth', dotted, domain= 1.15:4.8,variable=\t,smooth,samples=75]
plot ({-\t r}: {0.98*\t});		
\end{tikzpicture}
\end{center}
\caption{A 35-day cycle of \emph{wetonan} in the Javanese calendar. The regular Gregorian seven-day cycle is arranged clockwise from Monday to Sunday. The five-day \emph{pasaran} cycle is arranged from inside to outside. The white disk, red, blue, green, and yellow rings correspond to L\^{e}gi, Pahing, Pon, Wag\'{e}, and Kliwon, respectively.} \label{wetonanspiral}
\end{figure}

The \emph{wetonan} cycle combines the five-day \emph{pancawara} cycle with the seven-day \emph{wuku} week cycle. Each \emph{wetonan} cycle lasts for $7 \times 5 = 35$ days, with 35 distinct combinations of the couple `\emph{dinapitu} \emph{pasaran}' which is called the \emph{weton}. Figure~\ref{wetonanspiral} displays the 35-day \emph{wetonan} cycle of the dual \emph{dinapitu pasaran}. The seven-day \emph{wuku} cycle is arranged clockwise from Monday to Sunday and the five-day \emph{pasaran} day progresses from the center of the disk outwardly from Legi to Kliwon. It shows a spiral pattern that repeats seven times and sweeping five sectors each, indicated by solid black, solid blue, dashed-black, dashed-blue, dashed-dotted black, dashed-dotted blue, and dotted black spirals, respectively.

Although the \emph{weton} can be calculated independently using either equation~\eqref{weekgreg} or~\eqref{weekjav}, it is also possible to propose a single congruence formula:
\begin{equation}
\begin{aligned}
w &= k + \floor*{\frac{153(m + 1)}{5}} + 15Y + \floor*{\frac{Y}{4}}+ 19C + \floor*{\frac{C}{4}} + 5\\
W & \equiv w \, \tpmod{7} \\
P & \equiv w \, \tpmod{5}
\end{aligned}
\label{wetoneq}
\end{equation}

\noindent where $k$, $m$, $C$, and $Y$ are the same variables as defined for previous congruence formulas, $w$ is a 35-day offset congruence, and $W$ or $P$ corresponds to the index of \emph{dinapitu} or \emph{pasaran} day, respectively. Each term of the congruence can be analyzed as follows:
\begin{itemize}
\item $k$ represents the progression of the day of the week based on the day of the month, since each successive day results in an additional offset of one;
\item $+\floor*{\frac{153(m + 1)}{5}}$ adds the proper amount of days for each month, considering January and February as the 13\textsuperscript{th} and 14\textsuperscript{th} months of the previous year, respectively;
\item $+15Y$ adds an offset of $365 \tpmod{35} \equiv 15$ days for the common, non-leap years;
\item $+\left\lfloor {\frac{Y}{4}}\right\rfloor$ adds one more day for the leap years since $366 \tpmod{35} \equiv 16$;
\item $+19C$ adds $36,524\tpmod{35} \equiv 19$ days for any regular century (a century with non-leap year);
\item $+\floor*{\frac{C}{4}}$ adds one more day for a century leap year that is divisible by $400$ since $36,525\tpmod{35} \equiv 20$;
\item $+5$ adjusts the final offset to fit the indexes defined for \emph{dinapitu} $W$ and \emph{pasaran} $P$ after $\tpmod{7}$ and $\tpmod{5}$, respectively.
\end{itemize}

Let take the same example of 8 July 1633~CE with $k = 8$, $m = 7$, $C = 16$, and $Y = 33$. Using \eqref{wetoneq} we have $w \equiv 8 + 244 + 8 + 304 + 4 + 5 \equiv 573$, $W \equiv 573 \tpmod{7} \equiv 6$, and $P \equiv 573 \tpmod{5} \equiv 3$. Hence, the 8 July 1633~CE was a `Jemuwah Lêgi'. Notice that we are able to replace both the congruence formulas~\eqref{weekgreg} and~\eqref{weekjav} using our original single congruence relationship~\eqref{wetoneq}. Obtaining both \emph{dinapitu} and \emph{pasaran} days using a tabular method has been attempted by~\cite{arcienaga2020personal} and a \emph{Java} application based on the Indian calendar for calculating the Javanese calendar has been developed by~\cite{gislen2019calendars6}.

The \emph{wetonan} cycle is especially important for divinatory systems, celebrations, and rites of passage as birth or death. Commemorations and events are held on days considered to be auspicious. In particular, the \emph{weton} of birth is considered playing an important role in any individual personality, in a similar way as a zodiac sign does in Western astrology. The two \emph{weton} days of future spouses are supposed to determine their background nature compatibility and are used to compute the best date of marriage using a strict arithmetic formula~\citep{utami2019math}. It also figures in the timing of many ceremonies of \emph{slametan} ritual meal and many other traditional divinatory systems~\citep{utami2020ethnomathematics}. The eve of `Jumat Kliwon' is considered particularly popular and auspicious for magical and spiritual matters~\citep{darling2004marking,arcienaga2020personal}.

The anniversary of Javanese birthday occurs every 35 days, so about a thousand times in a century (exactly 1043 times). For a newborn baby, the first occurrence of \emph{weton}, i.e., aged 35 days, is named `selapanan' where the parents will cut hairs and nails of their child for the first time. For adults, the \emph{weton} of birth is considered as an eminent day, not festive but expressing humility and blissfulness. For example, a person may fast (sometimes also the day before and the day after), stop his commercial activity, avoid taking any major decision, or simply be more generous to surrounding people through philanthropic actions like share a blessed meal or some `jajan pasar' (early morning fresh sweet snacks from the market).

The \emph{weton} for the birth and death of Sultan Agung is `Jemuwah Lêgi', which is also the first day of the Javanese calendar he created. This \emph{weton} is therefore one of the recurrent noble days (see subsubsection \emph{Dina Mulya}) and  considered as an important night for pilgrimage. Indeed, the \emph{weton} of every King's birth is a special day; the present Sultan Hamengkubuwana X was born on 2 April 1946~CE, a `Selasa Wage'. This \emph{weton} has been chosen for His ascension to the throne, on 7 March 1989~CE. Every `Selasa Wage', the animated touristic center of Yogyakarta, Jalan Malioboro, is closed to motor vehicles, and the covered sidewalks `kaki lima', where commercial activity usually abound, is entirely cleaned.

As another especially prominent example, the present palace of Yogyakarta has been inhabited by Sultan Hamengkubuwana I and His regal suite on the `13 Sura 1682~AJ' (7 Oktober 1756~CE), a `Kemis Pahing'. In 2015, the Governor decided that every `Kemis Pahing', schoolchildren, public servants, and in particular those working in territorial services of Yogyakarta, must wear the traditional costume all day long, as a reminder of their regional culture. Merchants from traditional markets are encouraged to do the same.

\begin{table}
\caption{List of names of the 12 \emph{wulan} Javanese lunar months and the associated number of days, depending on the \emph{taun} and \emph{kurup} (see Table~\ref{taun}).}
\label{wulan}
\begin{tabular}{p{0.5cm}p{2.7cm}>{\centering\arraybackslash}p{1.8cm}>{\centering\arraybackslash}p{1.5cm}>{\centering\arraybackslash}p{1.5cm}>{\centering\arraybackslash}p{1.5cm}>{\centering\arraybackslash}p{1.5cm}}
\hline\noalign{\smallskip}
 &  & \multicolumn{5}{c}{Wulan length (days)}\\
 &  & {Taun 1-4,6-8} & \multicolumn{4}{c}{Taun 5}\\
No & Wulan Name & {All Kurup} & {Kurup 1} & {Kurup 2} & {Kurup 3} & {Kurup 4}\\
\noalign{\smallskip}\svhline\noalign{\smallskip}
1 & Sura         & \cellcolor{black!20}30 & \cellcolor{black!20}30 & \cellcolor{black!20}30 & \cellcolor{black!20}30 & \cellcolor{black!20}30\\
2 & Sapar       & 29 & 29 & \cellcolor{black!20}30 & \cellcolor{black!20}30 & 29\\
3 & Mulud        & \cellcolor{black!20}30 & \cellcolor{black!20}30 & 29 & 29 & \cellcolor{black!20}30\\
4 & Bakdamulud   & 29 & 29 & 29 & 29 & 29\\
5 & Jumadilawal  & \cellcolor{black!20}30 & \cellcolor{black!20}30 & \cellcolor{black!20}30 & 29 & \cellcolor{black!20}30\\
6 & Jumadilakir & 29 & 29 & 29 & 29 & 29\\
7 & Rejeb        & \cellcolor{black!20}30 & \cellcolor{black!20}30 & \cellcolor{black!20}30 & \cellcolor{black!20}30 & \cellcolor{black!20}30\\
8 & Ruwah        & 29 & 29 & 29 & 29 & 29\\
9 & Pasa         & \cellcolor{black!20}30 & \cellcolor{black!20}30 & \cellcolor{black!20}30 & \cellcolor{black!20}30 & \cellcolor{black!20}30\\
10 & Sawal        & 29 & 29 & 29 & 29 & 29\\
11 & Dulkangidah         & \cellcolor{black!20}30 & \cellcolor{black!20}30 & \cellcolor{black!20}30 & \cellcolor{black!20}30 & \cellcolor{black!20}30\\
12 & Besar     & 29/30 & \cellcolor{black!20}30 & \cellcolor{black!20}30 & \cellcolor{black!20}30 & 29\\
\noalign{\smallskip}\svhline\noalign{\smallskip}
 & {Total (days)} & 254/355 & 355 & 355 & 354 & 354\\
\noalign{\smallskip}\hline\noalign{\smallskip}
\end{tabular}
\end{table}

\subsubsection{Wulan}

The lunar month is named a \emph{wulan} and lasts for 29 or 30 days. There are 12 different names:  `Sura', `Sapar', `Mulud', `Bakdamulud', `Jumadilawal', `Jumadilakir', `Rejeb', `Ruwah', `Pasa', `Sawal', `Dulkangidah' (or `Sela'), and `Besar'. The length of each \emph{wulan} is attributed as follows (see Table~\ref{wulan} and the next subsubsections):
\begin{itemize}
\item `Sura', `Rejeb', `Pasa', and `Dulkangidah' are always 30 days;
\item `Bakdamulud', `Jumadilakir', `Ruwah', and `Sawal' are always 29 days;
\item `Sapar', `Mulud', `Jumadilawal', and `Besar' have lengths depending on the \emph{taun} and \emph{kurup}.
\end{itemize}

\subsubsection{Taun}

A \emph{taun} is a cycle on 12 \emph{wulan} and corresponds to the Javanese lunar year. There are eight different \emph{taun} names: `Alip', `Ehé', `Jimawal', `Jé', `Dal', `Bé', `Wawu', and `Jimakir', formed by different \emph{wulan} length sequences (see Table~\ref{wulan}). For seven of the \emph{taun}, the sequence alternates monotonically between 30 and 29-day lengths. For the final \emph{wulan} `Besar', it can be either 29 or 30 days, depending not only on the \emph{taun} but also on the \emph{kurup}, the 120 lunar year cycle, which has different day length sequences for the fifth \emph{taun} `Dal'. As a result, the total day length of a \emph{taun} varies from $354$ (short or normal year, named `Taun Wastu') to $355$ (long or leap year, named `Taun Wuntu'), as described in Table~\ref{taun}:
\begin{itemize}
\item `Alip', `Jimawal', and `Bé' are always normal years, with a \emph{wulan} `Besar' of 29 days;
\item `Jimawal' and `Wawu' are always leap years, with a \emph{wulan} `Besar' of 30 days;
\item `Jé' and `Dal' are normal or leap depending on the \emph{kurup};
\item `Jimakir' is a leap year for the 14 first \emph{windu}, but becomes a normal year for the final \emph{windu} of a \emph{kurup} cycle.
\end{itemize}
\begin{table}
\caption{List of names of the eight \emph{taun} Javanese lunar years forming a \emph{windu}, and the associated number of days, an alternate between 354 (short or normal year) and 355 (long or leap year, gray-shaded cells) days, depending on the \emph{kurup}.}
\label{taun}
\begin{tabular}{p{0.5cm}p{1.5cm}p{1.8cm}p{2.2cm}>{\centering\arraybackslash}p{1.2cm}>{\centering\arraybackslash}p{1.2cm}>{\centering\arraybackslash}p{1.2cm}>{\centering\arraybackslash}p{1.2cm}}
\hline\noalign{\smallskip}
 &  & &  & \multicolumn{4}{c}{Taun Length (days)}\\
No & Name & Krama & Meaning & {Kurup 1} & {Kurup 2} & {Kurup 3} & {Kurup 4}\\
\noalign{\smallskip}\svhline\noalign{\smallskip}
1 & Alip & Purwana & intention & 354 & 354 & 354 & 354\\
2 & Ehé & Karyana & action & \cellcolor{black!20}355 & \cellcolor{black!20}355 & \cellcolor{black!20}355 & \cellcolor{black!20}355\\
3 & Jimawal & Anama & work & 354 & 354 & 354 & 354\\
4 & Jé & Lalana & destiny & 354 & 354 & \cellcolor{black!20}355 & \cellcolor{black!20}355\\
5 & Dal & Ngawanga & life & \cellcolor{black!20}355 & \cellcolor{black!20}355 & 354 & 354\\
6 & Bé & Pawaka & back and forth & 354 & 354 & 354 & 354\\
7 & Wawu & Wasana & orientation & 354 & 354 & 354 & 354\\
8 & Jimakir & Swasana & empty & \cellcolor{black!20}355 & \cellcolor{black!20}355 & \cellcolor{black!20}355 & \cellcolor{black!20}355\\
\noalign{\smallskip}\hline\noalign{\smallskip}
\end{tabular}
\end{table}

Each \emph{taun} is assigned by a monotonic increasing number, based on the Indian calendar `Saka'. The reason was Sultan Agung decided to continue the counting from the Shalivahana era, which was 1555 at the time when inaugurating the Javanese calendar~\citep[cf.][]{nuraeni2017application}. Thus, the Javanese calendar began on `1 Sura Alip 1555~AJ', which corresponds to 8 July 1633~CE.

\subsubsection{Windu and Lambang}

Eight \emph{taun} make a \emph{windu}~\citep{proudfoot2006old}. Despite the variability of each \emph{taun} length, the total length of a normal \emph{windu} is constant since it always contains both five short and three long \emph{taun}, which is a total of 2,835 days (about 7 years 9 months in the Gregorian/Islamic calendar). This corresponds to exactly 81 \emph{wetonan}. This means that each `New Windu' day, dated as `1 Sura Alip', falls on the same \emph{weton}.

There is an exception to that rule: the final \emph{windu} of a \emph{kurup} cycle (see the next subsubsection) is always shortened by one day, with a 29 days \emph{wulan} `Besar' during the final \emph{taun} `Jimakir'. This induces a shift in the \emph{wetonan} cycle at each \emph{kurup}.

Furthermore, there are four different names of \emph{windu}: `Adi', `Kuntara', `Sengara', and `Sancaya' that compose a 32 \emph{taun} cycle. Another cycle of 16 \emph{taun} is combined using two different names of \emph{windu}, called \emph{lambang}: `Kulawu' and `Langkir'. These two cycles are summarized in Table~\ref{windu}.
\begin{table}
\caption{List of names of the four \emph{windu} and two \emph{lambang}.}
\label{windu}
\begin{tabular}{p{3.5cm}p{3.5cm}}
\hline\noalign{\smallskip}
Windu Name & Lambang Name\\
\noalign{\smallskip}\svhline\noalign{\smallskip}
Adi     & Langkir\\
Kuntara & Kulawu\\
Sêngara & Langkir\\
Sancaya & Kulawu\\
\noalign{\smallskip}\hline\noalign{\smallskip}
\end{tabular}
\end{table}

\subsubsection{Kurup}

The longest cycle in the Javanese calendar is called a \emph{kurup}, formed by 15 \emph{windu}, which is equivalent to 120 \emph{taun} or 1440 \emph{wulan}~\citep{gislen2019calendars4}. But the very last \emph{wulan} of the cycle, i.e., the twelfth \emph{wulan} `Besar' of the eighth \emph{taun} `Jimakir' of the fifteenth \emph{windu}, has only 29 days, such as the total length of a \emph{kurup} is $2,835 \times 15 - 1 = 42,524$ days (about 116 years and 6 months in the Gregorian/Islamic calendar)~\citep{rosalina2013aplikasi}. This is the same number of days as in 120 lunar years of the Tabular Islamic calendar. This is a rule-based variation of the Islamic Hijri calendar. Although the number of years and months are identical, the months are determined by arithmetical rules instead of observation or astronomical calculations.

Moreover, each \emph{kurup} determines:
\begin{itemize}
\item which of the \emph{taun}  `Jé' or `Dal' has a long \emph{wulan}  `Besar',
\item the sequence of \emph{wulan} lengths in the \emph{taun} `Dal',
\end{itemize}
\noindent as given in Tables~\ref{wulan} and~\ref{taun}. Hence, the full date sequences in the calendar vary between \emph{kurup}.

As the \emph{weton} of the first day of a \emph{kurup} repeats at each first day of the \emph{windu}, a \emph{kurup} is named using the corresponding \emph{weton} falling on `1 Sura Alip'. Table~\ref{kurup} lists the first five \emph{kurup}.

Meanwhile, the Sultanate of Mataram was divided under the Treaty of Giyanti between the Dutch and Prince Mangkubumi in 1755~CE. The agreement divided ostensible territorial control over Central Java between Yogyakarta and Surakarta Sultanates. The former was ruled by Prince Mangkubumi, also known as Raden Mas Sujana or Hamengkubuwono~I (1717--1792~CE) and the latter was administrated by Sinuhun Paliyan Negari, who was known as Pakubuwana~III (1732--1788~CE)~\citep{ricklefs1974jogjakarta,soekmono1981pengantar,frederick1993indonesia,brown2018short}.
\begin{table}
\caption{List of names of the five first \emph{kurup} Javanese 120 lunar year cycles, their short names (a contraction of the \emph{weton} at each new \emph{windu}, i.e., on `1 Sura Alip'), the first and last \emph{taun}, the total amount of \emph{taun}, and the starting dates in the Gregorian calendar.}
\label{kurup}
\begin{tabular}{p{0.4cm}p{1.8cm}p{1.5cm}>{\raggedleft\arraybackslash}p{2cm}>{\raggedleft\arraybackslash}p{2cm}>{\raggedleft\arraybackslash}p{0.8cm}>{\raggedleft\arraybackslash}p{2.4cm}}
\hline\noalign{\smallskip}
No &  Kurup Name & Short Name & First Taun (AJ) & Last Taun (AJ) & Taun & Start Date (CE)\\
\noalign{\smallskip}\svhline\noalign{\smallskip}
1 & Jamingiyah & A'ahgi & Alip 1555 & Jimakir 1674 & 120 & 8 July 1633\\
2 & Kamsiyah & Amiswon & Alip 1675 & Ehé 1748 & 74 & 11 December 1749\\
3 & Arbangiyah & Aboge & Jimawal 1749 & Jimakir 1866 & 118 & 28 September 1821\\
4 & Salasiyah & Asapon & Alip 1867 & Jimakir 1986 & 120 & 24 March 1936\\
5 & Isneniyah & Anenhing & Alip 1987 & Jimakir 2106 & 120 & 26 August 2052\\
\noalign{\smallskip}\hline\noalign{\smallskip}
\end{tabular}
\end{table}

During the second \emph{kurup} (1749--1821~CE), some experts realized that the Javanese calendar was still one day behind compared to the Islamic Hijri calendar. Hence, the King of Surakarta, Susuhunan Pakubuwana V (1784--1823~CE), decided to end the \emph{Kurup} `Amiswon' in the year 1748~AJ, even though it had only been running for nine \emph{windu} and two \emph{taun}. So, the \emph{taun} `Ehé' 1748~AJ, which was supposed to be a leap year, was made only 354 days and the third \emph{kurup} `Aboge' started on the \emph{taun} `Jimawal' 1749~AJ. But some noticed that it would be more appropriate if the incrementation of \emph{kurup} should have been carried out two lunar years before, namely on the \emph{taun} `Alip' 1747~AJ. As a consequence of this delay, the third \emph{kurup} `Aboge' is only 118 \emph{taun} long. However, the Sultanate of Yogyakarta did not make a similar decision and pursued the second \emph{kurup} normally, so that the calendar in the two concurrent regions was different during 46 years (see Table~\ref{kurupkesultanan}). On \emph{taun} `Jimakir' 1794~AJ, the Sultan of Yogyakarta, Hamengkubuwana VI (1821--1877~CE), finally agreed and decided that the third \emph{kurup} `Aboge' will also end with \emph{taun} `Jimakir' 1866~AJ, reconciling the two calendars.
\begin{table}
\caption{List of dates of the second and third \emph{kurup} in the Sultanate of Yogyakarta.}
\label{kurupkesultanan}
\begin{tabular}{p{0.4cm}p{1.8cm}p{1.5cm}>{\raggedleft\arraybackslash}p{2cm}>{\raggedleft\arraybackslash}p{2cm}>{\raggedleft\arraybackslash}p{0.8cm}>{\raggedleft\arraybackslash}p{2.4cm}}
\hline\noalign{\smallskip}
{\it No} & {\it Kurup Name} & {\it Short Name} & {\it First Taun (AJ)} & {\it Last Taun (AJ)} & {\it Taun} & {\it Start Date (CE)}\\
\noalign{\smallskip}\svhline\noalign{\smallskip}
2 & Kamsiyah & Amiswon & Alip 1675 & Jimakir 1794 & 120 & 11 December 1749\\
3 & Arbangiyah & Aboge & Alip 1795 & Jimakir 1866 & 72 & 16 May 1866\\
\noalign{\smallskip}\hline\noalign{\smallskip}
\end{tabular}
\end{table}

The fourth and present \emph{kurup} `Asapon' is planned to last a normal length of 120 lunar years, and will end on `29 Besar Jimakir 1986~AJ', which is 25 August 2052~CE. The following day will start the fifth \emph{kurup} `Anenhing' (`Alip Senen Pahin') on `1 Sura 1987~AJ Alip', but the sequence of leap years has not yet been decided, so it is formally impossible to calculate the exact \emph{dina}, \emph{wulan}, and \emph{taun} for such a distant future date. Nevertheless, there is no obstacle with other strictly monotonic cycles like \emph{weton}, \emph{wuku}, and \emph{windu}.

\subsubsection{Dina Mulya}

\emph{Dina mulya} are the noble days in the Javanese calendar. Excepted for the `Siji Sura' which is the new lunar year and falls on the first day of the first \emph{wulan} every \emph{taun}, others are associated with a specific \emph{weton} and specific \emph{taun} or \emph{wuku} (see Table~\ref{dinamulya}).
 
\begin{table}
\caption{List of the noble days \emph{dina mulya}.}
\label{dinamulya}
\begin{tabular}{p{2cm}p{2.2cm}p{2cm}p{0.7cm}p{0.7cm}p{0.7cm}p{2.5cm}}
\hline\noalign{\smallskip}
Dina Mulya & Weton & Wuku & Dina & Wulan & Taun & Occurrences\\
\noalign{\smallskip}\svhline\noalign{\smallskip}
Siji Sura & -- &  -- &  1 &  Sura &  -- & 1 every 354/355 days\\
                &   &   &   &      &   & (new lunar year)\\
\noalign{\smallskip}
Aboge & Rêbo Wage & -- & -- & -- & Alip & 10 during the {\it Taun}\\
                &   &   &   &      &   & every 7/8 years\\
\noalign{\smallskip}
Daltugi & Sêtu Legi & -- & -- & -- & Dal & 10 during the {\it Taun}\\
                &   &   &   &      &   & every 7/8 years\\
\noalign{\smallskip}
Kuningan & Sêtu Kliwon & Kuningan & -- & -- & -- & every 210 days\\
\noalign{\smallskip}
Hanggara Asih & Selasa Kliwon & Dukut & -- & -- & -- & every 210 days\\
\noalign{\smallskip}
Dina Mulya & Jemuwah Kliwon & Watugunung & -- & -- & -- & every 210 days\\
\noalign{\smallskip}
Dina Purnama & Jemuwah Lêgi & -- & -- & -- & -- & every 35 days\\
\noalign{\smallskip}\hline\noalign{\smallskip}
\end{tabular}
\end{table}

\subsection{Computer implementation of the Javanese calendar}

The computation of the full Javanese calendar has been implemented using \emph{GNU Octave} scientific language with a single function `weton.m' \citep{BeauducelWeton20}. When using a computer, however, the determination of the weekday or \emph{pasaran} day from a date in the Gregorian calendar does not require the congruence formula~\eqref{wetoneq}. In fact, most computer languages are able to calculate the exact number of days that last from any reference date, correctly taking into account leap years. In the case of the Javanese calendar, the linear timeline will be the number of days counted from 8 July 1633~CE, which falls on `Jemuwah Lêgi 1 Sura Alip 1555~AJ Jamingiyah', the first day of the first \emph{kurup}. Using that facility, calculation of the 7-day week cycle can be made by a simple modulo 7, the \emph{pasaran} cycle by a modulo 5, the \emph{wetonan} cycle by a modulo 35, and the \emph{pawukon} cycle by a modulo 210 functions. Moreover, since any modern digital calendar is able to set a periodic event, a specific \emph{weton} date repeated every 35 days will smoothly give the corresponding \emph{wetonan} cycle over the whole calendar.

On the other hand, the computation of \emph{dina}, \emph{wulan}, \emph{taun}, \emph{windu}, and \emph{kurup} is more complicated since the exact sequences vary throughout history following human decisions, making these cycles not exactly cyclic nor monotonic. Hence, the proposed computing strategy is to construct, for each \emph{kurup}, a \emph{windu} table as a $8 \times 12$ matrix of rows \emph{taun} versus columns \emph{wulan}, containing the day length of that specific \emph{wulan}. This matrix is repeated 15 times to form a complete \emph{kurup}, or less for the second and third \emph{kurup}. Then, all the matrices are concatenated. The cumulative sum of this table elements, in the row order, gives the total number of days that lasts from the origin at the beginning of each lunar month, and can be compared to the linear timeline described above. Thus, a simple `table lookup' function will give the corresponding \emph{kurup}, \emph{windu}, \emph{taun}, and \emph{wulan} indexes, and the \emph{dina} will be given by the remainder.

\section{Conclusion} 		\label{conclude}

In this chapter, we have discussed the cultural, historical, and arithmetic aspects of the Javanese calendar. Along with the internationally acknowledged Gregorian and the majority-embraced religiously Islamic Hijri calendars, the Javanese calendar has its unique place among the heart of many Javanese people in Indonesia as well as Javanese diaspora overseas. Although many Javanese people have adopted modern lifestyle, the Javanese calendar is still utilized in various daily affairs, including to choose the best possible time for arranging a wedding day. 

While determining the day of the week for any given date can be computed using Zeller's congruence algorithm of modulo seven, the \emph{pasaran} day of the Javanese calendar can be calculated using a new congruence formula of modulo five. Additionally, we have also proposed a unique and combined congruence formula for calculating both the day of the week and \emph{pasaran} day for any given date in the Gregorian calendar. Furthermore, using a computer program \emph{GNU Octave} `weton.m'~\citep{BeauducelWeton20}, all the cycles of the Javanese calendar, i.e., \emph{wetonan}, \emph{wuku}, \emph{wulan}, \emph{windu}, \emph{lambang}, and \emph{kurup} can also be determined straightforwardly.

\section*{Acknowledgment}
The authors would like to thank Matthew Arciniega (Vortx, Inc., Ashland, Oregon, United States of America) for sharing the contents of his old website and Roberto Rizzo (University of Milano--Bicocca, Italy) for pointing to the article written by~\cite{proudfoot2007search}. NK dedicated this chapter to his father, who introduced to and taught him about calendars during his early childhood. FB warmly thanks Alix Aimée Triyanti for her cultural influence and inspiration, and all his Javanese friends for their enthusiastic support.

\bibliographystyle{spbasic}  
\bibliography{2012author} 


\end{document}